\documentclass{article}
\title{Geometry of superficial elements}
\usepackage[german,english]{babel}

\usepackage{amsmath,amssymb,amsthm,amsfonts,latexsym}
\usepackage{amsthm}

\newtheorem{thm}{Theorem}[section]
\newtheorem{lem}[thm]{Lemma}
\newtheorem{prop}[thm]{Proposition}
\theoremstyle{definition}
\newtheorem{defi}[thm]{Definition}
\newtheorem{rem}[thm]{Remark}

\author{Romain Bondil\thanks{Address : Fakult\"at f\"ur Mathematik der Ruhr-Universit\"at, Universit\"atsstr. 150, Geb. NA 2/31, 44780 Bochum, Germany. Email : romain.bondil@ruhr-uni-bochum.de}}

\newcommand{\Z}{\mathbb{Z}}
\newcommand{\N}{\mathbb{N}}

\DeclareMathOperator{\proj}{Proj}
\DeclareMathOperator{\spec}{Spec}

\DeclareMathOperator{\Min}{Min}
\DeclareMathOperator{\Ann}{Ann}
\DeclareMathOperator{\Ass}{Ass} 

\DeclareMathOperator{\Supp}{Supp}

\DeclareMathOperator{\ini}{in}

\newcommand{\EI}{\mathcal{E}_I}

\newcommand{\OEI}{\mathcal{O}_{\mathcal{E}_I}}
\newcommand{\GIM}{\mathcal{G(M)}}
\newcommand{\fdEi}{(f)^\#_{\EI}}
\newcommand{\BIM}{\mathcal{B(M)}}

\newcommand{\Mc}{\mathcal{M}}

\newcommand{\Oc}{\mathcal{O}}

\newcommand{\OSI}{ \mathcal{O}_{S_I} } 

\begin{document}
\maketitle

\begin{abstract}
In this paper,\footnote{{\it 2000 M.S.C}.~: 13H15, 14C17
 } we study superficial elements of an ideal with respect to a module from a geometrical point of view, using blowing-ups. The notion of weak transform is particularly relevant to this study. We use this viewpoint to get a ``natural'' proof of a theorem by D. Kirby characterizing  those superficial elements. We also indicate how the same result may be algebraically derived from a more recent theorem of  Flenner and Vogel.
\end{abstract}
\section*{Introduction}

The notion of superficial element was introduced by P. Samuel in~\cite{Sa}  for the study of multiplicities of primary ideals in local rings.
 
As is well-known and recalled in \S~\ref{sec-Samuel-Kirby}, superficial elements are not only well-behaved with respect to multiplicities but also with respect to Hilbert polynomials.

More basically, superficial elements can be defined for {\em any} ideal simply by a property of {\em the kernel} of the multiplication by the class of this element in the associated graded ring (resp. module) and an  at first sight surprising fact (yielding the afore-mentioned nice behaviour w.r.t. numerical invariants) is that this property forces the {\em cokernel} of the same multiplication map to be very nice.

This property is part of a characterization of superficial elements given by D. Kirby, recalled in thm.~\ref{thm-Kirby}.

The purpose of this paper is to give a more picturesque approach to superficial elements using projective geometry and blowing-up schemes.

It is natural to translate the condition to be superficial as a property on the $\proj$ on the associated graded ring. However, one gets further understanding if one embeds this $\proj$ as the exceptional divisor of the blow-up scheme.

A first characterization of superficial elements on the blow-up is given at the end of \S~\ref{sec-first-geom-form} and is expressed in term of {\em the weak  transform} of these elements (cf. prop.~\ref{prop-first-char}).

In \S~\ref{sec-wk-vs-st}, we compare this weak transform to the more usual {\em strict} (or {\em proper} transform) giving a criterion of equality.

Section~\ref{sec-regularity} contains our  main result~: since superficiality is a condition of regularity on the exceptional divisor of the blow-up, i.e. on a Cartier subscheme, it may be used to define  a two-terms regular sequence, which, once  inverted,  give  two non-trivial regularity conditions, the second precisely connecting  the weak and strict transforms (thm.~\ref{thm-super-equiv-permu-reg}).

In Section~\ref{s-back-to-Kirby} we go back to Kirby's formulation showing that the two conditions found in our theorem~\ref{thm-super-equiv-permu-reg}· exactly give the two Kirby's conditions.

Eventually, we mention in the last section how Kirby's result may be also derived from a more general result by Flenner and Vogel which makes clear the connection between kernel and cokernel of certain maps between associated graded rings. However, the proof is not so direct in the case of superficial elements which may be zero divisors.

\smallskip


We end this introduction  by setting-up hypotheses and notation  valid throughout this paper:

\smallskip
\noindent {\bf Setup~--}Let $R$ always stand for  a noetherian ring, $I$ an ideal of $R$ and $M$ a finitely generated $R$-module. It is standard to consider the following graded rings and modules (cf. e.g. \cite{Ei} Chap.~5 for the ring and module structure)~:
\begin{eqnarray}
G(R):=G_I(R):=\underset{n\in \N}{\oplus}\:  I^n/I^{n+1},\\
\label{eq-def-GIM}
G(M):=G_I(M):=\underset{n\in \N}{\oplus}\: I^n M/I^{n+1}M,
\end{eqnarray}
associated to the filtration of $R$ by the $(I^n)_{n\in N}$.
We will also denote $G_n$ for the component of degree $n$ of a graded module $G$, and write  $G(s)$ for the graded module with ``shifted'' gradation   $[G(s)]_n:=G_{n+s}$.

\section{Superficial elements after Samuel and Kirby}
\label{sec-Samuel-Kirby}

\noindent {\bf Multiplication by $\bar{f}$~--} For an element $f\in I^s\setminus I^{s+1}$ one considers its so-called initial form $\bar{f}$ which is its class in $I^s/I^{s+1}$.

Multiplication by $\bar{f}$ defines a map of degree zero between graded modules~:
\begin{equation}
\label{eq-def-mfbar}
m_{\bar{f}}\: : \: G(M)\rightarrow G(M)(s),
\end{equation}

Being interested in the properties  of this map, one studies the associated exact sequence~:
\begin{equation}
\label{eq-exact-seq-mfbar}
0\rightarrow \ker m_{\bar{f}} \rightarrow G(M)\rightarrow G(M)(s)\rightarrow G(M)(s)/\bar{f}.G(M)\rightarrow 0.
\end{equation}

Considering also the graded module $G(M/fM)$  defined as in (\ref{eq-def-GIM}) replacing $M$ by $M/fM$, it is easy to check that, denoting by $\ini_I (fM)$ the graded submodule of $G(M)$ whose $n$-th component is~:
\begin{equation}
\label{eq-def-ini(fM)}
[\ini_I (fM)]_n:=(I^nM \cap fM)/I^{n+1}M,
\end{equation}
one has the isomorphism~:
$$G(M/fM)\cong G(M)/\ini_I (f.M).$$
In particular,  one has a natural surjection~:
\begin{equation}
\label{eq-surj-sur-G(I,fM)}
G(M)/\bar{f}G(M)\rightarrow G(M/fM),
\end{equation}
corresponding to the obvious inclusion $\bar{f}.G(M)\subset \ini_I(fM)$.\par

\medskip
\noindent{\bf Introduction of superficial elements~--}As he was interested in properties of Hilbert functions (see below), Samuel only considered properties valid for components of large enough degree and introduced (cf. \cite{Sa} p.~182, the following equivalent formulation is the one in \cite{AC}~VIII  \S~7 No.~5 p.~79)~:
\begin{defi}\label{def-super-bbki}
An element $f\in I^s\setminus I^{s+1}$  is said to be  {\em superficial} (of degree~$s$) for $M$ with respect to $I$ if, and only if, there is an integer $n_0$ such that the multiplication map $m_{\bar{f}}$  defined in~(\ref{eq-def-mfbar}) is  injective from  $G(M)_{n-s}\rightarrow G(M)_n$ for all $n\geq n_0$.

\end{defi}

A somewhat surprising property of superficial elements is the following~:
\begin{prop}
\label{prop-bon-ker-impliq-bonker}
Let $R$ be a noetherian ring and $M$ a finitely generated $R$-module.
Let $I$ be an ideal of $R$.
If $f\in I^s\setminus I^{s+1}$ is superficial for $M$ w.r.t. $I$ then the epimorphism~(\ref{eq-surj-sur-G(I,fM)}) is an isomorphism, and hence one has the following sequence~:
\begin{equation}
\label{eq-exact-seq-cas-super}
0\rightarrow G(M)(-s)\overset{m_{\bar{f}}}{\rightarrow} G(M)\rightarrow G(M/fM)\rightarrow 0,
\end{equation}
which is exact between homogeneous components of degree $n\geq n_0$.\footnote{the $n_0$ is intended to be the same as in def.~\ref{def-super-bbki}}
\end{prop}

From the view-point of comparison of graded properties of $M$ and $M/fM$, this miracle says that {\em injectivity of $m_{\bar{f}}$ is enough to get ``the right cokernel''}.

\medskip
\noindent {\bf Application to Hilbert functions~--}
If one is interested (as Samuel originally was) by properties of Hilbert functions, one  takes $R$ to be a noetherian local ring and $I$ an ideal such that the length  $l(M/IM)$ is finite, then one defines the Hilbert function of $G(M)$ as~:
$$H_{I,M}(n):=l(I^n M /I^{n+1}M)$$
and the exact sequence~(\ref{eq-exact-seq-cas-super}) gives the nice relationship for all $n\geq n_0$~:
$$H_{I,M/fM}(n)=H_{I,M}(n)-H_{I,M}(n-s),$$
which allows to compute this Hilbert function by induction on dimension.
 
\medskip
\noindent {\bf Characterization of superficial elements by prop.~\ref{prop-bon-ker-impliq-bonker}}

In fact, prop.~\ref{prop-bon-ker-impliq-bonker} is part of a full  characterization of superficial elements  by D. Kirby (in~\cite{Ki} thm.~3) which  we now state\footnote{In loc. cit. it was formulated only for $M=R$ but the generalization  to any $R$-module is straightforward.}. We use the following standard notation (cf. e.g. \cite{Ei} \S~3.6)~:  denote by $(0\: :_M\: I):=\{m\in M, mI=0\}$ and $\Gamma_I(M):=(0\: :_M \: I^\infty):=\cup_{n\in \N} (0\: :_M \: I^n)$.
 
\begin{thm}
\label{thm-Kirby}
Let $R$ be a noetherian ring, $I$ an ideal of $R$, and $M$ a finitely generated $R$-module.
An $f\in I^s\setminus I^{s+1}$ is superficial for $M$ w.r.t. $I$ if and only if both of the following conditions hold~:  

\noindent {\rm (i)} multiplication by $f$ from $M/\Gamma_I(M) \rightarrow M/\Gamma_I(M)$ is injective.

\noindent {\rm (ii)} the map in~(\ref{eq-surj-sur-G(I,fM)})~: $G(M)/\bar{f}G(M)\rightarrow G(M/fM)$ is an isomorphism between components of degree $n\geq n_0$ 
\end{thm}

\section{Superficial elements, definition on the blowup}
\label{sec-first-geom-form}
Before shifting to the language of projective geometry, we remind the reader about another piece of commutative algebra~:

\medskip
\noindent{\bf Associated primes~--} From the properties relating associated primes and non zero divisors (cf. \cite{Ei} thm.~3.1) we immediately get (cf. \cite{AC} VIII p.~79)~:
 
\begin{rem}
\label{rem-ass-elem-super}
With the same notation as in def.~\ref{def-super-bbki} and considering the set $\Ass(G(M))$ of homogeneous prime ideals in  $G(R)$ associated to $G(M)$ (cf. \cite{Ei} Chap.~3) the condition {\em $f$ is superficial}  is equivalent to the condition that {\em $\bar{f}$ does not belong to the $p_i\in \Ass(G(M))$ such that $p_i$ does not contain $G(R)_1=I/I^2$}.
\end{rem}


\medskip
\noindent {\bf Projective Formulation~--} Both def.~\ref{def-super-bbki} and rem.~\ref{rem-ass-elem-super} are more naturally formulated using projective geometry.

Denote  $\EI:=\proj(G(R))$, then the $G(R)$-graded module $G(M)$ defines a (sheaf of) $\OEI$-module $\GIM:=\widetilde{G(M)}$ using the standard functor $\; \widetilde{ }\;$ (cf. e.g. \cite{Ii} \S~3.4), which also gives  a map~:
\begin{equation}
\label{eq-mfbar-tilde}
\widetilde{m_{\bar{f}}}\: : \: \GIM  \rightarrow  \GIM(s)
\end{equation}
associated to the map  $m_{\bar{f}}$ defined in~(\ref{eq-def-mfbar}).

A basic fact  about this functor $\; \widetilde{ }\;$ is (cf. e.g. \cite{Ii} \S~7.1 or \cite{EGA} II \S~2.7)~:

\begin{lem}
\label{lem-TN-equiv-tilde}
If $R$ is a noetherian ring and $M$ is a finitely generated module, then the injectivity of the map $m_{\bar{f}} : G(M)_n \rightarrow G(M)_{n+s}$ for $n$ large enough  is equivalent to the injectivity of $\widetilde{m_{\bar{f}}}\: : \: \GIM \rightarrow \GIM(s)$.
\end{lem}




In the same spirit, the subset $\Ass(\GIM)$ of points in $\EI$  associated to the $\OEI$-module $\GIM$ (cf.~\cite{EGA} IV 3.1) precisely corresponds to  the elements of $\Ass(G(M))$ fulfilling the condition in rem.~\ref{rem-ass-elem-super} (called {\em projectively relevant}).

Finally, consider the subscheme of $\EI$ defined by the homogeneous ideal $\bar{f}. G(R)$ that we denote by $(f)^\#_{\EI}$ (this notation is to be explained below). With  the foregoing remarks, we reformulate def.~\ref{def-super-bbki} and rem.~\ref{rem-ass-elem-super} as follows~:

\begin{lem}
\label{lem-proj-super}
Let $R$ be a noetherian ring, $I$ an ideal of $R$ and $M$ a finitely generated $R$-module. Let $f\in I^s\setminus I^{s+1}$. Then 
$f$ is superficial for $M$ w.r.t. $I$  if, and only if, (two equivalent formulations)

\noindent  {\rm (i)} $\widetilde{m_{\bar{f}}}\: : \: \GIM \rightarrow \GIM(s)$ is injective, 


\noindent {\rm (ii)} $\Supp (f)^\#_{\EI}\cap \Ass(\GIM)=\emptyset$, 
where $\Supp$ denotes the underlying set of the scheme $(f)^\#_{\EI}$.
\end{lem}

\medskip
\noindent {\bf Embedding $\EI$ on the blowup $S_I$~--}
To get a better {\em geometric} understanding of both the map $\widetilde{m_{\bar{f}}}$ and the scheme $(f)^\#_{\EI}$ in lem.~\ref{lem-proj-super}, we may  embed the scheme~$\EI$ in the blowup scheme $S_I:=\proj B(R)$, where~:
$$B(R):=R[It]=\underset{n\in \N}{\oplus} I^n t^n,$$
\noindent is graded by the powers of $t$, with the convention $I^0=R$.

For a subscheme of a projective scheme, it is possible to consider either its global ``homogeneous'' equations (i.e. the homogeneous ideal defining it)  or its local ones in affine charts. As far as affine charts are concerned, fixing a basis~$(h_0,\dots,h_r)$ of the ideal $I$ of $R$, one defines an affine open covering of $S_I$ by the $U_i:=D_+(h_i t)\cong \spec R[I/h_i]$.

Denote also $b_I : S_I \rightarrow \spec R$ the blowup morphism induced by the inclusion of $R$ in $B(R)$.

Then $\EI$ is counter-image $(b_I)^{-1}(V(I))$ and hence globally defined by the homogeneous ideal $I.B(R)$ in  $S_I$ and locally by the equation $h_i$ in each  $U_i$ for~$i=0,\dots,r$.

\medskip
\noindent {\bf Total and weak transforms~--}
Taking an $f\in I$, one defines its {\em total transform} $(b_I)^{*}(f)$ as the counter-image $(b_I)^{-1}(V(f))$ hence globally defined by the homogeneous principal ideal $f.B(R)$ and locally also by $f$ in each chart $U_i$.


Before going further, we put the emphasis on the following piece of terminology that will be of some importance later~:
\begin{defi}
\label{defi-loc-princ-cartier}
i) If $X$ is a scheme and $Y$ is a subscheme of $X$, we say that $Y$ is {\em locally principal} if there is a covering of $X$ by affine open subsets $U_i=\spec A_i$ such that $Y\cap U_i$ is defined by a principal ideal $(f_i)\subset A_i$.

ii) We say that $Y$ is a {\em Cartier subscheme} if it fulfills condition i) and further the $f_i\in A_i$ are all non zero divisors.
\end{defi}

Making no special assumption on $f\in R$, we see that $(b_I)^*(f)$ is simply a locally principal subscheme of $S_I$ whereas $\EI$ is a Cartier subscheme.
This is enough for the following~:

\begin{defi}
\label{defi-weak-transf}
We define the {\em weak transform} $(f)^\#$ of an element $f\in I^s \setminus I^{s+1}$ on the blowup $S_I$ of $I$ in $\spec R$, as the subscheme~:
$$(f)^\#:=(b_I)^*(f)-s\,  \EI,$$
\noindent where the $-$ sign means that one takes locally the quotient of the equations in each affine chart.\footnote{\label{ft-weak-transf} The terminology {\em weak transform} may be found in \cite{Hk} p.~142. I chose the  $^\#$ sign to mean that this weak transform can be bigger than the {\em strict} transform of $(f)$.}

Remark that this weak transform $(f)^\#$ is also globally defined by the homogeneous ideal $f t^s B(R)$ of $B(R)$.
\end{defi}


\medskip 
\noindent {\bf The scheme $\fdEi$ is the pull-back of the weak transform~--}
Considering the homogeneous equation of $(f)^\#$ and the morphism (of degree zero) of graded rings $B(R) \rightarrow G(R)=B(R)/IB(R)$, that sends $f t^s$ on $\bar{f} \in I^s/I^{s+1}$, one may define the intersection scheme $(f)^\#\cap \EI$\footnote{also called the pull-back on $\EI$ of the locally principal subscheme $(f)^\#$}  as the subscheme of $\EI$ defined by the homogeneous ideal $\bar{f}G(R)$.

Hence, $\EI\cap (f)^\#$ is nothing but the scheme $\fdEi$ of lem.~\ref{lem-proj-super}, so that (ii) of this lemma may be re-phrased as in the following~:

\begin{prop}
\label{prop-first-char}
Under the same notation as in def.~\ref{def-super-bbki}, $f\in I^s\setminus I^{s+1}$ is superficial for $M$ with respect to $I$ if, and only if, on the blowup $S_I$ of $\spec R$ along $I$, denoting $\mathcal{G(M)}=\widetilde{G(M)}$, and $(f)^\#$ the weak transform defined in~\ref{defi-weak-transf}, we have~:
\begin{equation}
\label{eq-def-super-geom}
\Supp\: ((f)^\#)\cap \Ass\: (\mathcal{G(M)})=\emptyset.
\end{equation}
where $\Supp$ denotes the locus on $S_I$ defined by the weak transform   $(f)^\#$ and $\Ass$ the (finite) set of associated points of $\GIM$ (included in $\EI\subset S_I$).

In the special case $M=R$, (\ref{eq-def-super-geom}) reads~:
\begin{equation}
\label{eq-cas-M-egal-R}
\Supp\: ((f)^\#)\cap \Ass\: (\EI)=\emptyset.
\end{equation}
\end{prop}

\section{Weak transform vs. strict transform}
\label{sec-wk-vs-st}

In def.~\ref{defi-weak-transf}, we introduced the {\em weak transform} $(f)^\#$ of an element $f\in I^s\setminus I^{s+1}$ on the blowup $S_I$ of the ideal $I$  in $\spec R$.

Let us recall the classical definition of strict (or proper) transform of a subscheme (cf. e.g. \cite{E-H} p.~168). Recall first that the scheme-theoretic closure of a subscheme is by definition the smallest closed subscheme containing it.

\begin{defi}
\label{defi-strict-transf}
Let $Y\subset \spec R$ a closed subscheme and $b_I\: : S_I \rightarrow S$ the blowup of $I$ in $S=\spec  R$. The strict transform $Y'$ of $Y$ by $b_I$ is by definition the scheme-theoretic closure of the counter image of $Y\setminus V(I)$ by $b_I$, which we denote by~:
$$Y':=\overline{(b_I)^{-1}(Y\setminus V(I))}^{sch}.$$
Because of the uniqueness of the structure of subscheme on an open subset, this subscheme is also the scheme-theoretic closure of $(b_I)^{-1}(Y)\setminus (b_I)^{-1}(V(I))=(b_I)^{-1}(Y)\setminus \EI$. For the same reason the scheme-structure of $V(I)$ is not to be taken into account in this definition, but simply its support.
\end{defi}

In the particular case of $Y=(f):=\spec R/(f)$, its strict transform is~:
$$(f)':=\overline{(f)^*\setminus \EI}^{sch},$$
but since,  with $(f)^\#$ the weak transform of def.~\ref{defi-weak-transf}, we trivially have $(f)^*\setminus \EI=(f)^\#\setminus \EI$, one may just as well say~:
$$(f)'=\overline{(f)^\#\setminus \EI}^{sch}.$$

This gives in particular the inclusion of schemes $(f)'\subset (f)^\#$ (whence the $^\#$ notation cf. footnote~\ref{ft-weak-transf}), and the condition of equality is given by the following elementary lemma~:
\begin{lem}
\label{lem-weak-equal-strict}
With the same notation as in \S~\ref{sec-first-geom-form} on the blowup, the weak transform  $(f)^\#$ is equal to the strict transform $(f)'$ if, and only if, the support of the exceptional divisor $\EI$ does not contain any associated point to the scheme $(f)^\#$,  what we denote by~:
\begin{equation}
\label{eq-weak-equal-strict}
\Supp (\EI)\cap \Ass (f)^\#=\emptyset.
\end{equation}
\end{lem}
\begin{proof}
This is a simple application of the theory of primary decomposition~:
the inclusion  of closed subschemes $(f)'\subset (f)^\#$ is proper if, and only if, one has a decomposition~:
$$(f)^\#=(f)'\cup X_I,$$ 
into closed subschemes, where $X_I$ contains a point associated to $(f)^\#$ (non trivial decomposition), which must be in the support of $\EI$ since $(f)'$ and $(f)^\#$ coincide outside $\EI$.
\end{proof}

Because of the similarity of conditions in prop.~\ref{prop-first-char} (\ref{eq-cas-M-egal-R}) and in lem.~\ref{lem-weak-equal-strict}, up to the permutation of the roles of $\Supp$ and $\Ass$, we will investigate the precise connection between these two conditions in the next section (see the conclusive remark~\ref{rem-cas-M-R}).
\smallskip

\noindent {\bf Case of modules~--} Let $M$ be a (finitely generated)  $R$-module and denote by $\BIM$ the $\OSI$-module defined on the blowup $S_I$ of $I$ in $\spec R$ by the graded module~:
$$B(M)=\underset{n\in \N}{\oplus}\, I^n M t^n.$$
 
Then, for $f\in I$ one may define the ``restrictions'' of $\BIM$ to the strict transforms $(f)'$ and weak transform $(f)^\#$ by~:
\begin{equation}
\label{eq-def-st-transf-BIM}
\BIM_{(f)'}:=\BIM \otimes_{\OSI} \Oc_{(f)'},
\end{equation}

\noindent (resp. $\BIM_{(f)^\#}$ by tensor product with $\Oc_{(f)^\#}$).

In this context, the more general form of lemma~\ref{lem-weak-equal-strict} (with the same proof, i.e. primary decomposition for modules) is~:
\begin{lem}
\label{lem-module-wk-trsf-egal-st}
The two $\OSI$-modules $\BIM_{(f)'}$ and $\BIM_{(f)^\#}$ defined in~(\ref{eq-def-st-transf-BIM}) are isomorphic if, and only if,
$$\Supp(\EI)\cap \Ass(\BIM_{(f)^\#})=\emptyset.$$
\end{lem}

\section{Superficiality and regularity}
\label{sec-regularity}

We first recall  standard constructions in  projective geometry (cf. \cite{EGA} II \S~2.6 or \cite{Ii} \S~7.1)~:

\noindent {\bf Morphism $\alpha$~--} Let $A:=\oplus_{n\in \N}  A_n$ be a graded ring and $M=\oplus_{n \in \Z} M_n$ a graded $A$-module. Then one defines for all $n\in \Z$ a morphism (of $A_0$-modules)~:
\begin{equation}
\label{eq-def-alpha}
\alpha_n \: : \: M_n \rightarrow \Gamma(\proj A, \tilde{M}(n)),
\end{equation}
by, locally in each chart $D_+(f_i)$ of $\proj A$, sending $m\in M_n$ simply to~$m/1 \in \Gamma (D_+(f_i),\tilde{M}(n))$.
\medskip

\noindent {\bf Multiplication by sections of $\Oc_X(s)$~--}
Take again $A$ to be a graded ring, $X=\proj A$ and $\sigma \in \Gamma(X,\Oc_{X}(s))$. Such a global section defines for any  $\Oc_X$-module $\Mc$, a morphism of multiplication by $\sigma$~:
\begin{equation}
\label{eq-mult-par-section}
m_\sigma: \Mc \rightarrow \Mc(s).
\end{equation}
\medskip
Explicitly, assume for simplicity that $A=A_0[A_1]$ so that charts $D_+(f_i)$ with $f_i \in A_1$ cover $X$. Then, in each $U_i:=D_+(f_i)$, $\sigma =a_i f_i^s$ with $a_i\in \Oc_X(U_i)$.

\begin{rem}
\label{rem-def-regular}
In particular, since multiplication by $f_i^s$ is certainly injective from $\Mc(U_i)\rightarrow \Mc(s)(U_i)$  the injectivity of $m_\sigma$ is equivalent to the elements $a_i \in \Oc_X(U_i)$ being $\Mc(U_i)$ regular, with the usual terminology of  a regular element for a module (cf. e.g. \cite{Ei} Chap.~17).
If this is the case we will say that the section $\sigma\in \Gamma(X,\Oc_X(s))$ is $\Mc$-regular.
\end{rem}

\medskip
\noindent{\bf Application to $X=\proj(G(R))$~--} Now we return to the setting of lemma~\ref{lem-proj-super}~: for $G(R)=\oplus_{n\in \N} I^n/I^{n+1},\:  \EI:=\proj G(R)$ and $\bar{f}\in G(R)_s$ one gets from (\ref{eq-def-alpha}) a global section $\alpha_s(\bar{f})\in \Gamma(\OEI, \OEI(s))$ and it is direct from the definitions that the morphism $\widetilde{m_{\bar{f}}} : \GIM \rightarrow \GIM(s)$ considered in the cited lemma, coincides with the multiplication by $\alpha_s(\bar{f})$ in the sense of (\ref{eq-mult-par-section}).

With the terminology of rem.~\ref{rem-def-regular}, we get the following avatar of lemma~\ref{lem-proj-super}~:
\begin{lem}
\label{lem-f-super-equiv-alpha-gim-regular}
In the same setting  as in lemma~\ref{lem-proj-super} and using  the foregoing definitions~:
$f$ is superficial for $M$ w.r.t. $I$ if, and only if, the section $\alpha_s(\bar{f})$ is $\GIM$-regular.
\end{lem}

\medskip
\noindent {\bf Application to $X=S_I=\proj B(R)$~--} As we did in \S~\ref{sec-first-geom-form}, we now shift from $\EI$ to the larger space $S_I=\proj(\oplus _{n\in \N}\, I^n t^n)$ of blowup of $I$ in $\spec R$.

 The morphism $\alpha$ defined in (\ref{eq-def-alpha}) above also applies to get a global section $\alpha(ft^s) \in  \Gamma(S_I, \OSI(s))$. Recall that we denote $\BIM$ for the $\OSI$-module associated to $B(M):=\oplus_{n\in \N}\,  I^n M t^n$.  

Using all this, our main result boils down to the following permutation of elements in a regular sequence~:

\begin{thm}
\label{thm-super-equiv-permu-reg}
Let $R$ be a noetherian ring, $I$ an ideal of $R$ with basis $(h_0,\dots,h_r)$.
Denote $U_i=D_+(h_it)$ the corresponding affine open subsets in the blowup space $S_I=\proj(B(R))$.

Let $M$ be a finitely generated $R$-module and $\BIM$ the corresponding $\OSI$-module as defined above.

Then for an element $f\in I^s \setminus I^{s+1}$, the following four conditions are equivalent~:

\noindent {\rm a)} $f$ is superficial for $M$ w.r.t. $I$ (cf. lem.~\ref{lem-f-super-equiv-alpha-gim-regular}, or lem.~\ref{lem-proj-super}), 

\noindent {\rm b)} $\forall\, i=0,\dots,r,\:\forall x\in S_I, \: (h_i,f/h_i^s)$ is a $\BIM_x$-regular sequence (where ``$_x$'' denotes the localization at $x$),

\noindent {\rm c)} $\forall\, i=0,\dots,r,\:\forall x\in S_I, \: (f/h_i^s,h_i)$ is a $\BIM_x$-regular sequence,

\noindent {\rm d)} the following two conditions are satisfied ~:
 
{\rm (i)} $\alpha_s(ft^s)$ is $\BIM$-regular (definition  in rem.~\ref{rem-def-regular}),

{\rm (ii)} $\Supp(\EI)\cap\Ass(\BIM_{(f)^\#})=\emptyset$, where $(f)^\#$ is the weak transform introduced in~\ref{defi-strict-transf} and we consider the ``restriction'' modules in the sense of~(\ref{eq-def-st-transf-BIM}).
\end{thm}

\begin{proof}

Let us fix an $i\in \{0,\dots,r\}$, and work in the corresponding open subset $U_i\cong  \spec R[I/h_i]$ on the blowup.

\noindent $\bullet$ (a)$\Leftrightarrow$(b)~:

Since  $\BIM(U_i)=M[I/h_i]$, $h_i$ is  regular for any $R$-module $M$, and hence (b) reduces to the condition~: $f/h_i^s$ is regular for $\BIM(U_i)/h_i\BIM(U_i)=\GIM(U_i)$ at all $x\in U_i$.  But, considering the corresponding  classes $\bar{f}$ and $\bar{h_i}$ in $G(R)_s$ and $G(R)_1$ respectively, this is the same as 
$\bar{f}/\bar{h_i}^s$ being regular for $\GIM(U_i)$.

We now recognize the local equation of $\alpha_s(\bar{f})$ and hence the condition $\alpha_s(f)$ is $\GIM$-regular of lemma~\ref{lem-f-super-equiv-alpha-gim-regular}, whence (b)$\Leftrightarrow$(a).

\noindent $\bullet$ (b)$\Leftrightarrow$(c)~:
standard permutation property for regular sequences over a local ring, cf. e.g. \cite{Ei} cor.~17.2.

\noindent $\bullet$ (c)$\Leftrightarrow$(d)~: in the sequence $(f/h_i^s,h_i)$ in $U_i$ the first element defines the equation of $\alpha_s(ft^s)$ whence (i) in d) and the condition for the $h_i$ to be regular with respect to the quotient $\BIM(U_i)$ by $f/h_i^s$ translates as (ii) in d) since
$(f)^\#$ is the subscheme defined by  $\alpha_s(ft^s)$ and using always the same relation between associated points and zero divisors.\end{proof}

In the case $M=R$ we get the better sounding formulation~: 
\begin{rem}
\label{rem-cas-M-R}
In the special case $M=R$ in the foregoing theorem we get that 

\noindent (a) $f$ is superficial w.r.t. $I$

if, and only if,

\noindent (d) (i) the weak transform $(f)^\#$ is a Cartier subscheme of $S_I$ (cf. def.~\ref{defi-loc-princ-cartier}), and 
 
(ii) $\Supp(\EI)\cap \Ass((f)^\#)=\emptyset$, which is exactly the condition~(\ref{eq-weak-equal-strict}) in  lem.~\ref{lem-weak-equal-strict} so that the weak transform equals the strict transform.
\end{rem}

\section{Back to Kirby's theorem}
\label{s-back-to-Kirby}

We claim now that thm.~\ref{thm-Kirby} is easily recovered from our thm.~\ref{thm-super-equiv-permu-reg} condition (d). Precisely conditions (i) (resp. (ii)) correspond  in both results, as we now check~:

\begin{lem}
The following two conditions are equivalent~:

\noindent $(*)$  Condition (i) in thm.~\ref{thm-Kirby} i.e. multiplication by $f$~: $ M/\Gamma_I(M) \rightarrow M/\Gamma_I(M)$ is injective,

\noindent $(**)$ Condition (i) in thm.~\ref{thm-super-equiv-permu-reg}~(d) i.e.  $\alpha_s(ft^s)$ is $\BIM$-regular on the blowup scheme $S_I$.
\end{lem}
\begin{proof}
Condition $(**)$ is equivalent to the multiplication by $ft^s$~: $I^n M t^n \rightarrow I^{n+s} M t^{n+s}$ being injective for $n$ large, and since $t$ is certainly a non-zero divisor, we may just as well consider multiplication by $f$~: $I^n M \rightarrow I^n M$.

Hence to prove the lemma, one has to check for $n$ large, the following~:
\begin{equation}
\label{eq-egal-Ass}
\Ass_R (I^n M)=\Ass (M/\Gamma_I(M)).
\end{equation}

To prove (\ref{eq-egal-Ass}), one may reduce to the case $\Gamma_I(M)=0$.
Indeed, by noetherian condition $\Gamma_I(M)=(0:_M I^{n_1})$ for a certain $n_1$ and hence $I^{n_1}.M=I^{n_1}.(M/\Gamma_I(M))$.
So replacing $M$ by $M/\Gamma_I(M)$ and taking $n\geq n_1$, we may assume $\Gamma_I(M)=0$.

We now prove $\Ass(I^n M)=\Ass(M)$ for all $n$  in the case $\Gamma_I(M)=0$.

From the inclusion $I^n M \subset M$, one always has $\Ass(I^n M) \subset \Ass(M)$.

Conversely, if $p\in \Ass(M)$, one may localize at $p$ and denote $M$ and $R$ for $M_p$ and $R_p$. Then $p=\Ann (m)$ is equivalent to $p.m=0$, since $p$ is the maximal ideal.
 
Then, as $\Gamma_I(M)=0$ we have for any fixed $n$, $I^n.m\neq 0$ i.e. there is an $i\in I^n$ such that $i.x\neq 0$.

Hence, $i.p.x=0$ gives $p=\Ann(ix)$ i.e. $p\in \Ass(I^n.M)$.\end{proof}

\begin{lem}
The following two conditions are equivalent~:

\noindent $(*)$ Condition {\rm (ii)} in Kirby's thm.~\ref{thm-Kirby} i.e. the epimorphism~: 
$G(M)/\bar{f}.G(M) \rightarrow G(M/f.M)$ is an isomorphism between components of large degree,

\noindent $(**)$ Condition {\rm (ii)} in our thm.~\ref{thm-super-equiv-permu-reg}~(d) i.e. $\Ass(\BIM_{(f)^\#})\cap \Supp(\EI)=\emptyset$.
\end{lem}
\begin{proof}
From lemma~\ref{lem-module-wk-trsf-egal-st}, we know that condition $(**)$ is exactly the condition for the isomorphism of the two $\OSI$-modules defined by the restriction $\BIM_{(f)^\#}$ and $\BIM_{(f)'}$ of $\BIM$ to the weak and strict transform respectively.

But since these two modules are clearly isomorphic at each point outside $\EI$ it is equivalent to check that their pull-back on $\EI$ are actually isomorphic.
These pull-back (taking $\otimes \EI$) are exactly $\widetilde {G(M)/\bar{f}.G(M)}$ for the weak transform and $\widetilde{G(M/fM)}$ for the strict transform.

This is equivalent to condition $(*)$ from the standard lemma~\ref{lem-TN-equiv-tilde}.
\end{proof}

\section{More reasons for the miracle÷: a theorem of Flenner-Vogel}

The miracle referred to in the title is the one mentioned after prop.~\ref{prop-bon-ker-impliq-bonker}.

One may also get Kirby's characterization (thm.~\ref{thm-Kirby}) from the following more precise result in~\cite{F-V}; recall first that for an $R$-module $M$ one defines the {\em cycle} $Z(M)$  associated to $M$ by~:
$$Z(M)=\sum_{p\in \Min(M)} l(M_p)[p],$$
where $\Min(M)$ is the set of minimal associated primes to $M$.\par

\begin{thm}
\label{thm-F-V}
Let $R$ be a noetherian ring, $I$ an ideal of $R$ and 
$$0\rightarrow M_0 \rightarrow M_1 \rightarrow M_2 \rightarrow 0,$$
an exact sequence of finitely generated $R$-modules.

Considering the associated graded modules, one gets a complex $G(M_0)\rightarrow G(M_1) \rightarrow G(M_2)$ in which the last map is still onto and one has the equality of cycles~:
$$Z(\ker(G(M_0) \rightarrow G(M_1)))= Z(\ker (G(M_1)/G(M_0) \rightarrow G(M_2))).$$\end{thm}

Considering now an $f\in I^s\setminus I^{s+1}$, multiplication by $f$ in $M$  gives an exact sequence~:
$$0\rightarrow \Ann_M(f) \rightarrow M \rightarrow M \rightarrow M/fM \rightarrow 0,$$
that one has to split in order to apply thm.~\ref{thm-F-V} to the  two sequences~:
\begin{eqnarray}
\label{eq-FV1}
0 \rightarrow \Ann_M (f) \rightarrow M \rightarrow M/\Ann_M(f) \rightarrow 0,\\\label{eq-FV2} 
0 \rightarrow M/\Ann_M(f) \rightarrow M \rightarrow M/f.M \rightarrow 0.
\end{eqnarray}

The interested reader may deduce Kirby's theorem from theorem~\ref{thm-F-V} in the following (perhaps not so geometrically telling) way~:

\noindent (i) the conclusion in Kirby's theorem says that $\ker\big( G_I(\Ann f) \rightarrow G_I(M)\big)$ and $\ker\big( G_I(M)/G_I(fM)\rightarrow G_I(M/fM)\big)$ are irrelevant,

\noindent (ii) one has the following (short!) exact sequence
\begin{eqnarray*}
0\rightarrow \ker\big( G_I(M)/G_I(\Ann(f)) \rightarrow G_I(fM)\big) \rightarrow \\
\ker\big( G_IM/G_I(\Ann f) \rightarrow G_I(fM)\rightarrow G_I M\big)\rightarrow \ker\big( G_I(fM)\rightarrow G_IM \big)\rightarrow 0,
\end{eqnarray*}

\noindent (iii) The condition for $f$ to be superficial is that the middle-term of the sequence in (ii)  is irrelevant, which is equivalent to the two others being irrelevant. Then,  by theorem~\ref{thm-F-V} applied to (\ref{eq-FV1}) and (\ref{eq-FV2}), if these two terms are irrelevant one has the two conditions of  Kirby's theorem~\ref{thm-Kirby} because of (i).

\smallskip

{\footnotesize{\noindent{\bf Acknowledgment}~--The paper is a souped-up version of the first chapter of the author's Ph.D. thesis under L\^e D\~ung Tr\'ang (cf. \cite{Bo}).  The motivation for such a geometric study was mainly to complete the point of view in \cite{B-L}, where in spite of the title, we focused not on superficial elements but on so-called $v$-superficial elements of ideals (developed in the following chapters of the thesis). I would like to thank Mark Spivakovsky for his remarks as a referee for my  dissertation  and H. Flenner for his suggestions relative to the last parts of this work.}}

\end{document}